\documentclass[letterpaper, 10 pt, conference]{ieeeconf}  

\IEEEoverridecommandlockouts                              

\usepackage{cite}
\usepackage{amsmath,amssymb,amsfonts}
\usepackage{amsmath,bm}
\usepackage{algorithmic}
\usepackage{graphicx}
\usepackage{textcomp}
\usepackage{varioref}
\usepackage{wrapfig}
\usepackage{threeparttable}
\usepackage{dcolumn}
\newcolumntype{d}{D{.}{.}{-1}}
\makeglossary
\usepackage{fancyvrb}
\fvset{fontsize=\footnotesize,xleftmargin=2em}
\usepackage{lettrine}
\usepackage[colorlinks]{hyperref}
\usepackage{adjustbox}
\usepackage{graphicx}
\usepackage{epstopdf}
\usepackage{attachfile}
\usepackage{anyfontsize}
\usepackage{booktabs}
\usepackage{array}
\usepackage{varwidth}
\usepackage[section]{placeins}
\usepackage[intoc]{nomencl}
\usepackage{ifthen}
\usepackage{relsize}
\usepackage{cleveref}
\graphicspath{{./Figures/}} %

\title{\LARGE \bf Optimal Control of Wave Energy Converters Using Epsilon-Trig Regularization Method*
}

\author{Kshitij Mall$^{1}$ and Ehsan Taheri$^{1}$
\thanks{*This work was not supported by any organization.}
\thanks{$^{1}$Kshitij Mall and Ehsan Taheri are with the Department of Aerospace Engineering, Auburn University,
        Auburn, AL 36849, USA
        {\tt\small \{mall,etaheri\}@auburn.edu}}%
}

\begin{document}

\maketitle
\thispagestyle{empty}
\pagestyle{empty}
\graphicspath{{./Figures/}}

\begin{abstract}
The wave energy converter (WEC) devices provide access to a renewable energy source. Developing control strategies to harvest maximum wave energy requires solving a constrained optimal control problem. It is shown that singular control arcs may constitute part (or the entire) of extremal trajectories. Characterizing the optimal control structure, especially with the possibility of many switches between regular and singular control arcs, is challenging due to lack of \textit{a priori} information about: 1) optimal sequence as well as number of the regular and singular control arcs, and 2) the corresponding optimal switch times (from a regular to a singular arc and vice versa). This investigation demonstrates the application of a recently developed construct, the Epsilon-Trig Regularization Method (ETRM), to the problem of maximizing energy harvesting for a point-absorber WEC model in the presence of control constraints. Utility of the ETRM for the WEC problem is demonstrated by comparing its high-quality results against those in the literature for a number of test cases. 
\end{abstract}
\section{INTRODUCTION}
We are considering the problem of harvesting maximum wave energy from devices that make use of the surface motion of the waves \cite{babarit2006optimal}. The most commonly known devices in this category are \textit{point}- and \textit{linear}-absorber models. Control strategies to maximize energy absorption of WECs are achieved by solving a constrained optimal control problem (OCP), which is an active area of research \cite{hals2011constrained,scruggs2013optimal,li2014model,abdelkhalik2016control,coe2017comparison,wang2018review}. OCPs are traditionally solved using \textit{direct} and \textit{indirect} \cite{betts1998survey,silva2010smooth,mall2017epsilon,taheri2018generic} methods. 

The maximum wave energy harvesting problem of a point-absorber WEC has a control-affine Hamiltonian structure. For such Hamiltonian systems, a frequent phenomenon that \textit{may} occur is the appearance of \textit{singular} arcs, which usually complicates the solution procedure. The coefficient of the control input in affine-control Hamiltonian systems is called the \textit{switching function}. The sign of the switching function may alternate between positive and negative values thereby leading to the switching of the control input. Singular arcs, however, correspond to the cases in which the switching function vanishes for finite time intervals. In such cases, the Pontryagin's Minimum Principle (PMP) is not sufficient to characterize the extremal control and additional steps have to be taken. However, the mere existence of a control-affine structure in the Hamiltonian does not necessarily mean that the optimal control will consist of singular arcs \cite{stengel2012optimal}. 

A common practice for solving problems with mix regular-singular control structure involves a number of steps including 1) utilization of higher-order optimality conditions to find the control associated with a singular arc \cite{betts2010practical} and 2) dividing the entire problem into optimally connected segments of known regular and singular arcs. The algebraic expression for a singular control is obtained by taking successive time derivatives of the switching function until the control appears explicitly, which is a tedious task and is problem dependent. Additional conditions have also to be satisfied to ensure that the resulting singular control minimizes the Hamiltonian \cite{kelley1964second}. 

The aforementioned challenges can be overcome through the Epsilon-Trig Regularization Method (ETRM) \cite{mall2017epsilon}. The key step in this method is to alter the equations of motion (EOMs) of the OCP using trigonometric functions. This modification leads to significant consequences such that both regular and singular control arcs can be realized in a straightforward fashion. Another appealing feature is achieved by adopting the ETRM, namely, the original multi-point boundary-value problem (MPBVP) is reduced to a two-point boundary value problem (TPBVP); the latter is remarkably easier to solve. 

The main contribution of this work is the application of the ETRM to the WEC problem, which is known to have extremal solutions that consist of regular and singular control arcs. However, ETRM makes the numerical solution significantly easier such that a standard boundary-value solver such as MATLAB's \textit{bvp4c} can be used to solve these challenging problems. 

The remainder of this paper is organized as follows. Section \ref{sec:wecmodel} describes the point-absorber WEC model used in this study. A discussion of the TPBVP formulation and solution process for the WEC problem is given in Section \ref{sec:trigonomerizationmc}. Section \ref{sec:results} demonstrates the results and provides a comparison of solutions with those in the literature. Finally, concluding remarks are given in Section \ref{sec:conc}. 
\section{POINT ABSORBER MODEL FOR WAVE ENERGY CONVERTER}
\label{sec:wecmodel}
\begin{figure}[!ht]
\centering
\includegraphics[width=2.2in]{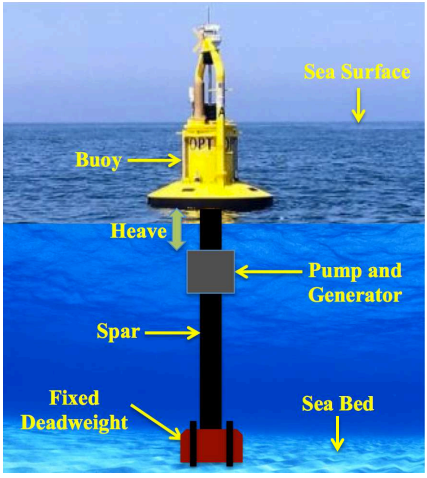}
\caption{A schematic of the WEC point absorber model.}
\label{fig:pointabsorber}
\end{figure}

Figure \ref{fig:pointabsorber} depicts the schematic for a typical point-absorber WEC model, where hydraulic cylinders are attached to a buoy \cite{tutorials2018alternative, falnes2002ocean}. The motion of the waves creates a vertical motion in the buoy, which results in pushing the hydraulic cylinders. These cylinders then drive the hydraulic motors, which in turn drive a generator. The power take-off (PTO) systems comprising of the hydraulic cylinders and motors thus translate the oscillating motion of the buoy to useful electrical energy. 

\begin{figure}[!ht]
\centering
\includegraphics[width=2.5in]{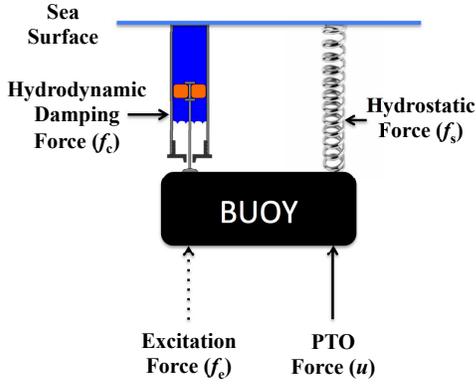}
\caption{Dynamic forces in the WEC point absorber model.}
\label{fig:wecdyn}
\end{figure}

Figure \ref{fig:wecdyn} depicts the various forces on the buoy. The dynamics of the WEC involve four forces: 1) a hydrostatic force, $f_{s}$, 2) a hydrodynamic damping force, $f_{c}$, 3) an excitation force, $f_{e}$, and 4) a PTO force.  The following assumptions are used in this study: 1) a linear dynamic model corresponding only to the heave motion is used, 2) the PTO force is chosen as the control force, which is assumed to act in the opposite direction to the heave motion, and 3) the wave frequency dependence of the hydrodynamic damping force is neglected.

Let $x$ denote the vertical displacement of the buoy. For a buoy with mass $m$, the dynamical model can be expressed as 
\begin{align} 
     m\ddot{x} &= f_{s} + f_{c} + f_{e} - u,\label{eqn:weceqm}
\end{align}
where the hydrostatic force, $f_{s} = -kx$, is the buoyancy force on the buoy, which is similar to the spring force with a spring constant $k$. The hydrodynamic force, $f_{c} = -c\dot{x}$, acts similar to the damping force on the system with a hydrodynamic damping constant, $c$. The excitation force, $f_{e} = \sum_{i=1}^{n} A_i\sin(\omega_i t + \phi_i)$, is the pressure effect around the immersed buoy or the float system (a periodic excitation force is considered here). In the relation for $f_e$, $A_{i}$ and $\phi_{i}$ are the amplitude and the phase for the frequency $\omega_{i}$, respectively, $t$ is the time, and $n$ is the total number of terms. 

The objective is to maximize the extracted energy, $E$, over a prescribed time interval, $t \in [t_0.t_f]$, which can be written as
\begin{align}
  E &= \int_{t_0}^{t_f} u(t)\dot{x}(t)\text{d}t.\label{eqn:objorig}
\end{align}
\section{WAVE ENERGY CONVERTER PROBLEM FORMULATION AND SOLUTION PROCESS}
\label{sec:trigonomerizationmc}
Since we employ PMP, we re-write the cost functional from Eq.~\eqref{eqn:objorig} to the form given in Eq.~\eqref{eqn:wectrigobj}, which is subjected to EOMs as shown in Eq.~\eqref{eqn:wectrigeoms}. The resulting OCP is given as
\begin{subequations}
\label{eqn:wec}
    \begin{align}
  J &= -\int_{t_{0}}^{t_{f}}u(t)x_{2}(t)\text{d}t,\label{eqn:wectrigobj}\\
  \dot{x}_{1} &= x_{2},~\dot{x}_{2} = \dfrac{f_{\text{e}} + f_s +f_c - u}{m},~\dot{x}_{3} = 1,\label{eqn:wectrigeoms}
\end{align}
  \end{subequations}
where $x_1$ is the displacement, $x_2$ is the velocity magnitude, $x_3$ is the time, and $\gamma$ is the control magnitude such that $u \leq |\gamma|$. In Eq.~\eqref{eqn:wectrigeoms}, the excitation force, $f_e$, can be approximated as a periodic or a non-periodic function. The ETRM regularizes this OCP using two simple trigonometric modifications given as
\begin{subequations}
    \begin{align}
u &=  \gamma\sin u_{\text{TRIG}},\label{eqn:wectrigu}\\
\dot{x}_1 &= x_2 + \epsilon\cos u_{\text{TRIG}},\label{eqn:wectrigeomnew}
\end{align}
  \end{subequations}
where $\epsilon$ is the error parameter used to regularize $u$ and $u_{\text{TRIG}}$ is the new control. The Hamiltonian associated with the regularized system can be written as
\begin{align} \label{eqn:wectrigham}
H &= - \gamma x_{2} \sin u\textsubscript{TRIG} + \lambda_{x_{1}}(x_{2} + \epsilon\cos u\textsubscript{TRIG})\\
& \ \ + \dfrac{\lambda_{x_{2}}(f_{e} - kx_{1} - cx_{2}  - \gamma\sin u\textsubscript{TRIG})}{m}  + \lambda_{x_{3}}. \nonumber 
\end{align}
Using the first order necessary conditions of optimality, also known as the Euler-Lagrange equations, the EOMs for the costates can be formed as shown in Eq.~\eqref{eqn:wectrigcostates}. The EOM for the costate $\lambda_{x_{3}}$ depends on the derivative of $f_{e}$ with respect to $x_{3}$. The periodic form of $f_{e}$ for cases 1 and 2 with $n = 5$, and the non-periodic form of $f_{e}$ for case 3 with $n = 8$ are discussed in the next section.
\begin{subequations}
\label{eqn:wectrigcostates}
\begin{align}
    \dot{\lambda}_{x_{1}} &= \dfrac{k\lambda_{x_{2}}}{m}, \label{eqn:wectrigcostate1}\\
\dot{\lambda}_{x_2} &= -\lambda_{x_{1}} + \dfrac{c\lambda_{x_{2}}}{m} + \gamma\sin u\textsubscript{TRIG}, \label{eqn:wectrigcostate2}\\
\dot{\lambda}_{x_{3}} &= -\frac{\lambda_{x_{2}}\sum_{i=1}^{n} A_{i}\omega_{i}\cos(\omega_{i}x_{3} + \phi_{i})}{m}. \label{eqn:wectrigcostate3}
\end{align}
\end{subequations} 

The switching function for this problem, $H_{1}$, is shown in Eq.~\eqref{eqn:wecswitch}. Note that $H_1$ is the switching function associated with the control in the original (non-regularized) problem. The optimal control law is given in Eq.~\eqref{eqn:wectrigcontrol} using the Euler-Lagrange equations, which is dependent on $H_{1}$. Even if the value of $H_{1}$ vanishes, the optimal control can be explicitly found from among these two control options using the PMP.  
\begin{subequations}
\label{eqn:wectrigcontrolswitch}
\begin{align}
H_{1} &= 
\dfrac{-(\lambda_{x_{2}} + mx_{2})}{m},\label{eqn:wecswitch}\\
u_{\text{TRIG}}^{*} &= 
\begin{cases}
\arctan\left(\dfrac{\gamma H_{1}}{\epsilon\lambda_{x_{1}}}\right),\\\\
\arctan\left(\dfrac{\gamma H_{1}}{\epsilon\lambda_{x_{1}}}\right) + \pi.\label{eqn:wectrigcontrol}\\
\end{cases}
\end{align}
\end{subequations}

\section{Results}
\label{sec:results}
In order to demonstrate the utility of the ETRM, three problem cases are shown and discussed in this study. The differences are mainly due to the type of boundary conditions enforced on the initial position and velocity, bounds on the control, and the type of excitation force. These cases are selected to represent a range of possible scenarios for the WEC problem and allow us to compare the results with those reported in Ref. \cite{zou2017optimal}. For the WEC problem $m$ is 2$\times$10\textsuperscript{5} kg, $k$ is 1.2$\times$10\textsuperscript{5} kg/s\textsuperscript{2}, and $c$ is 10\textsuperscript{5} kgm\textsuperscript{2}/s\textsuperscript{3}. Table \ref{table:weccases} summarizes the three cases for the WEC problem. 

\begin{table}[!ht]
 \begin{center}
           \caption{Summary of the considered cases.}
 \begin{tabular}{ccccc}
 \hline
 \hline
        \textbf{Case } & \textbf{Initial }   & \textbf{Control}  & \textbf{Excitation }\\
          \#               & \textbf{Conditions} &  \textbf{Bound} ($\times$10\textsuperscript{5} N)& \textbf{Force Type} \\
        \hline
         1                 & Fixed & $\mid\gamma\mid$  $\leq$ 1.5 & Periodic \\
         
         2                 & Free  & $\mid\gamma\mid$  $\leq$ 1.5 & Periodic\\
         
         3                 & Free & $\mid\gamma\mid$  $\leq$ 1.0 & Non-Periodic\\
         \hline
\end{tabular}
          \label{table:weccases}
 \end{center}
\end{table}

A comparison is drawn between this study and Ref. \cite{zou2017optimal} for cases 1 and 2. Since case 3 uses the fitted non-periodic function, no comparison is made with Ref. \cite{zou2017optimal}, which uses non-periodic data from real experiments. In addition, brute-force application of methods that use variants of direct optimization methods fails for problems with singular control arcs unless remedial actions are taken. In order to clarify this point, the solution obtained using a pseudo-spectral method (PSM) is given for case 1. 

In the numerical simulations corresponding to cases 1 and 2, $f$\textsubscript{e\textsubscript{P}} is the periodic excitation force as shown in Eq.~\eqref{eqn:wectrig5} with an amplitude vector, $A$\textsubscript{P}, described in Eq.~\eqref{eqn:wectrig6}. The frequency vector, $\omega$\textsubscript{P}, based on a periodic time period, $T_{P}$, is shown in Eq.~\eqref{eqn:wectrig7}. The value of $T_{P}$ is equal to 10 s, which is consistent with the literature \cite{zou2017optimal}. The phase vector of the excitation force, $\phi$\textsubscript{P}, is shown in Eq.~\eqref{eqn:wectrig8}.
\begin{subequations}
  \begin{align}
      f\textsubscript{e\textsubscript{P}} &= \sum_{i=1}^{5} A\textsubscript{P}\textsubscript{i}\sin(\omega\textsubscript{P}\textsubscript{i}x_{3} + \phi\textsubscript{P}\textsubscript{i}),\label{eqn:wectrig5}\\
      A\textsubscript{P} &= [1, 0.1, 0.03, 0.5, 0.01] \times 10^{5} ~\text{(N)}, \label{eqn:wectrig6}\\
\omega\textsubscript{P} &= \left[\dfrac{2\pi}{T\textsubscript{P}}, \dfrac{0.5\pi}{T\textsubscript{P}},  \dfrac{12\pi}{T\textsubscript{P}},  \dfrac{4\pi}{T\textsubscript{P}},  \dfrac{0.1\pi}{T\textsubscript{P}}\right] \text{(rad/s)},\label{eqn:wectrig7}\\
\phi\textsubscript{P} &= \left[\dfrac{\pi}{2}, \dfrac{\pi}{8}, \dfrac{\pi}{5}, \dfrac{\pi}{3}, \dfrac{\pi}{4}\right] \text{(rad)}.\label{eqn:wectrig8}
\end{align}
\end{subequations}    
Equation \eqref{eqn:fenew1} describes the trigonometric fit of the non-periodic excitation force, $f$\textsubscript{e\textsubscript{NP}}, used in case 3 of this study, where the value of $A$\textsubscript{NP} is 4$\times$10\textsuperscript{5} N. The constant vectors, $a$\textsubscript{NP}, $\omega$\textsubscript{NP}, and $\phi$\textsubscript{NP} are used to derive $f$\textsubscript{e\textsubscript{NP}} and are shown in Eqs. (\ref{eqn:fenew2})--(\ref{eqn:fenew4}), respectively. Since the OCT relies on derivatives of the state EOMs, discrete data for $f$\textsubscript{e\textsubscript{NP}} cannot be used directly to solve the WEC problem. Therefore, a continuous fitting (Fourier approximation) function as shown in Eq.~\eqref{eqn:fenew1} is used instead. 
\begin{subequations} 
 \label{eqn:fenew}
   \begin{align}
 f\textsubscript{e\textsubscript{NP}} = &A\textsubscript{NP}\sum_{i=1}^{8}a\textsubscript{NP}\textsubscript{i}\sin(\omega\textsubscript{NP}\textsubscript{i}x_{3} + \phi\textsubscript{NP}\textsubscript{i}),\label{eqn:fenew1}\\
a\textsubscript{NP} = &[6.255, 24.1, 0.4027, 1.511,\notag\\
&\ 0.3596, 0.9603, 0.6938, 20.71]~\text{(m)}, \label{eqn:fenew2}\\
\omega\textsubscript{NP} = &[0.6837, 0.7458, 1.354, 0.5228,\notag\\
&\ 1.054, 0.3953, 0.3246, 0.7512]~\text{(rad/s)},\label{eqn:fenew3}\\
\phi\textsubscript{NP} = &[0.4082, 1.727, -0.4019, -1.737,\notag\\
&\ -2.663, -1.51, -2.364, 4.73]~\text{(rad)}.\label{eqn:fenew4}
\end{align}
\end{subequations} 
The subsequent subsections include the results obtained using the ETRM for the three cases of the WEC problem (see Table \ref{table:weccases}). Structures of the optimal controls are also characterized in terms of a sequence of their underlying control arcs, i.e., `B' and `S' shorthand notations are used to represent bang and singular control arcs, respectively.


We have adopted a numerical continuation process\cite{grant2014rapid,taheri2016enhanced,junkins2018exploration} with two continuation sets for the three cases. Using this continuation approach, a simpler OCP with a time duration of $t_f = 1$ second and a higher value for the error parameter, $\epsilon$, is solved initially. In other words, a two-parameter family of OCPs are formed and the problem is solved by using a standard homotopy method. Note that one of the homotopy parameters, $t_f$, is a natural boundary condition on the problem, whereas $\epsilon$ is the control regularization parameter. In the first continuation set, this simpler initial solution serves as an initial guess for a subsequent complex problem comprising a longer time duration. The first continuation set is completed when the terminal time condition specified in the original problem is reached after a specified number of steps. The subsequent continuation set operates on reducing $\epsilon$ to a reasonably small value. The computation times shown for the three cases in Table~\ref{table:wecsummary} include the times required to complete these continuation sets.  

\subsection{Numerical Results for Case \#1}
Table~\ref{table:wectrigcase1} lists the boundary conditions for this case. According to Ref. \cite{zou2017optimal}, for the considered dynamics and in the presence of a periodic excitation, the initial conditions for the displacement, $x_{1,0}$, and the velocity, $x_{2,0}$, can be written as
\begin{align} 
x_{1,0}= -\frac{1}{2c}\sum_{i=1}^{5} \dfrac{A\textsubscript{i}\cos(\phi\textsubscript{i})}{\omega\textsubscript{i}},
x_{2,0} = \frac{1}{2c}\sum_{i=1}^{5} A\textsubscript{i}\sin(\phi\textsubscript{i}).
\end{align}

\begin{figure}[!htbp]
\centering
\includegraphics[width=3in]{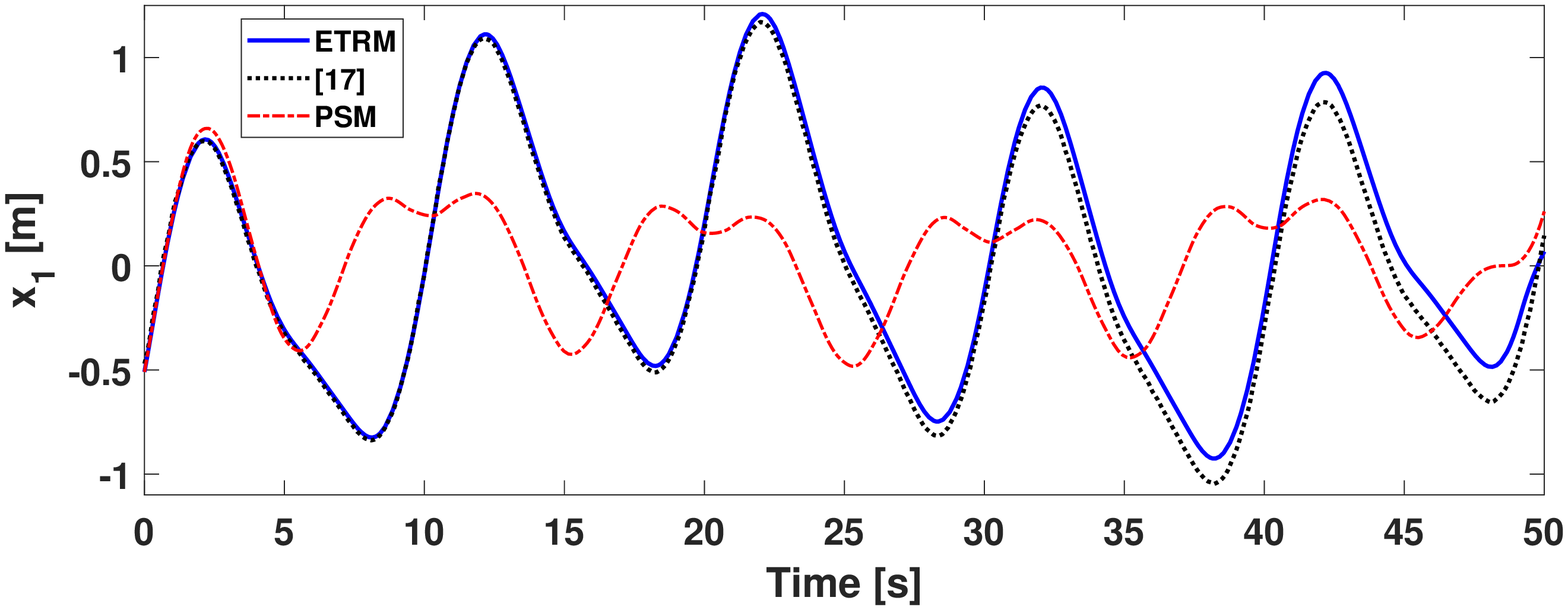}
\caption{Time history of $x_1$ for Case 1.}
\label{fig:case1_statex1}
\end{figure}

\begin{figure}[!htbp]
\centering
\includegraphics[width=3in]{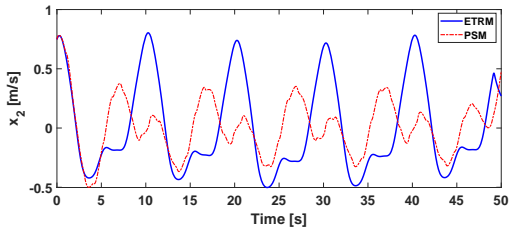}
\caption{Comparison of time histories of $x_2$ for Case 1.}
\label{fig:case1_statex2}
\end{figure}

\begin{table}[!ht]
 \begin{center}
           \caption{Initial and final conditions for case 1.}
 \begin{tabular}{p{2.5cm}cc}
 \hline
 \hline
        \textbf{Attribute} & \textbf{Initial Value} & \textbf{Final Value}\\\hline
        \textbf{Time (s)} & 0 & 50\\
        \textbf{Displacement (m)} & -0.5093 & Free\\
        \textbf{Velocity (m/s)} & 0.7480 & Free\\\hline
\end{tabular}
          \label{table:wectrigcase1}
 \end{center}
\end{table}
Figure \ref{fig:case1_statex1} shows $x_{1}$ obtained using the ETRM, Ref. \cite{zou2017optimal} and the PSM; Fig.~\ref{fig:case1_statex2} shows $x_{2}$ obtained using the ETRM and the PSM. The solutions for the ETRM and Ref. \cite{zou2017optimal} match well for the most part. However, the state solutions for the PSM are found to be completely spurious. 

\begin{figure}[!ht]
\centering
\includegraphics[width=3.4in]{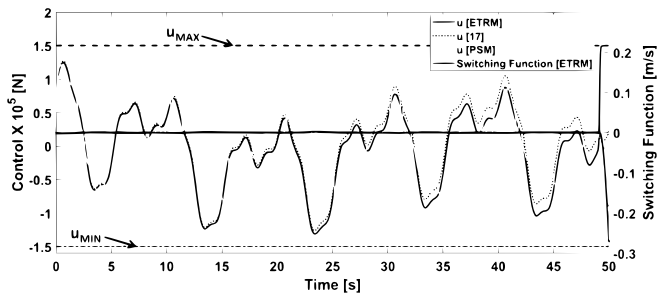}
\caption{Comparison of the time history of control inputs and the switching function between different methods for Case 1.}
\label{fig:case1_control}
\end{figure}

The control solutions for this case using the ETRM, Ref. \cite{zou2017optimal} and the PSM are shown in Fig.~\ref{fig:case1_control}. The control solution resulting from the ETRM is singular for the entire trajectory except for a small part in the end, where it assumes a bang form and attains the maximum value. Thus, the optimal control associated with the ETRM has an S-B sequence. Traditionally, one would have to solve a 3-point BVP using the OCT, which is more complicated than solving this case using the ETRM.
On the other hand, control solution from Ref. \cite{zou2017optimal} is purely singular, which is not the optimal strategy. In fact, the optimal control profile of the ETRM with a final bang harvests 0.841 MJ of energy, which is 5\% greater than the harvested energy in Ref. \cite{zou2017optimal} (0.8 MJ). The PSM has a jittery control solution, which is expected as the PSM has issues in solving singular control problems. When proper strategy is adopted, the PSM would be able to converge to the optimal solution. However, our goal is to show that this problem is indeed challenging. The switching function profile obtained using the ETRM (plotted in Fig.~\ref{fig:case1_control}) confirms the observations regarding the optimal control profile. Initially, the switching function stays at near-0-values corresponding to the singular control and then gains negative values corresponding to the bang control.  

\begin{figure}[!ht]
\centering
\includegraphics[width=3.5in]{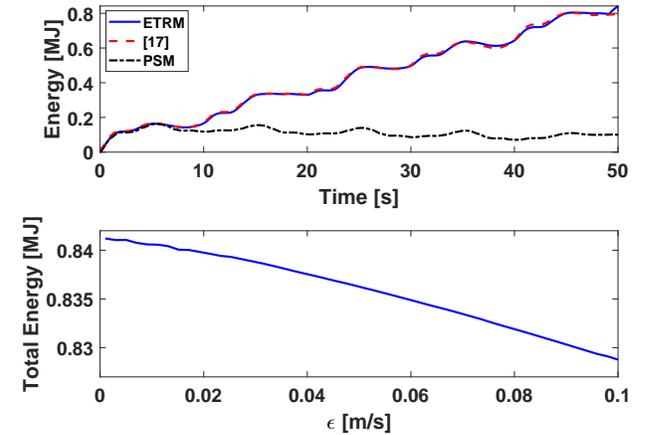}
\caption{Comparison of the time histories of the harvested energy (upper plot) and harvested energy vs. $\epsilon$ plot for Case 1.}
\label{fig:case1_energy}
\end{figure}

Figure \ref{fig:case1_energy} shows the energy time history plots obtained using the ETRM and from Ref. \cite{zou2017optimal}. When the value of $\epsilon$ is decreased from an initial chosen value of 0.1 m/s to 0.003 m/s using numerical continuation for the ETRM, higher and more accurate values are obtained for the harvested wave energy. We further note that the energy obtained using Ref. \cite{zou2017optimal} has higher values than the energy obtained using the ETRM at certain points along the trajectory (e.g., $t =41$ s). Since the objective is to increase the total energy over the entire time interval and the optimal control law using the ETRM involves S-B structure, the costates are different at certain points as compared to Ref. \cite{zou2017optimal}. In fact, the larger the control bound, the greater the differences would become since the optimal control input can take larger values.   
\subsection{Numerical Results for Case \#2}
\begin{figure}[!ht]
\centering
\includegraphics[width=3in]{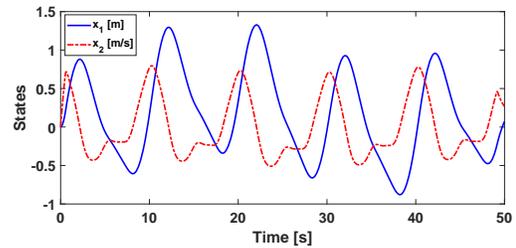}
\caption{Time histories of the states for Case 2.}
\label{fig:case2_states}
\end{figure}

\begin{figure}[!ht]
\centering
\includegraphics[width=3in]{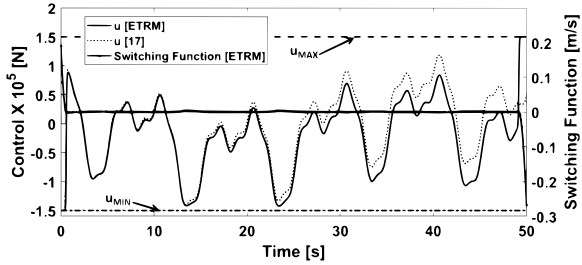}
\caption{Comparison of the time histories of the optimal control and switching function associated with the ETRM for Case 2.}
\label{fig:case2_control}
\end{figure}

\begin{figure}[!ht]
\centering
\includegraphics[width=2.8in]{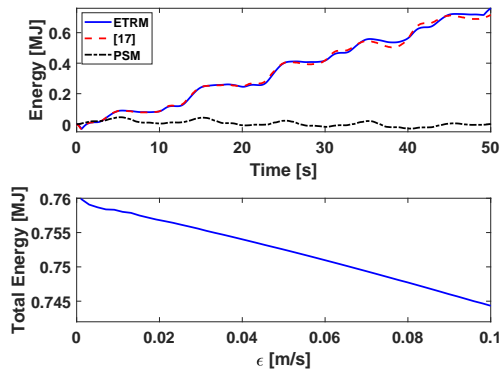}
\caption{Comparison of the time histories of harvested energies and changes in the harvested energy vs. $\epsilon$ for Case 2.}
\label{fig:case2_energy}
\end{figure}

For this case the initial values of $t$, $x_{1}$, and $x_{2}$ are 0 each and the final value of $t$ is 50 s. Figure \ref{fig:case2_states} shows the time history of the states for the ETRM and Ref. \cite{zou2017optimal}.
As shown in Fig.~\ref{fig:case2_control}, the control solution obtained using the ETRM starts with a bang form with a minimum value, switches to a singular form, and finally attains the bang form with a maximum value. Thus, the optimal sequence of control is B-S-B. The switching function matches with this observation as it initially has positive values corresponding to the first bang control. It then stays at 0 value corresponding to the singular control. Finally, it attains negative values corresponding to the second bang control. The early part of the optimal control profile obtained using the ETRM is identical to the control profile of Ref. \cite{zou2017optimal}, however, the final bang arc leads to a better energy harvest. Traditionally, one would have to solve a 4-point BVP using the OCT, which is more complicated than case 1. 
Figure \ref{fig:case2_energy} shows the energy time history plot for this case. The value of $\epsilon$ for the ETRM is decreased in an exact manner as case 1, resulting in higher and more accurate values of the harvested wave energy. The harvested energy associated with the ETRM control profile for this case is 0.76 MJ, whereas the energy harvested by the control profile of Ref. \cite{zou2017optimal} is 0.71 MJ. Thus, a 6.5\% improvement is achieved by using the ETRM due to the additional final bang arc. This seemingly negligible improvement gains considerable importance in WEC arrays or ``farms'', where such an improvement becomes multi-fold for a large number of WECs. Note that the method proposed in Ref. \cite{zou2017optimal} is indeed a simple strategy that may lead to sub-optimal control strategies, but it attains a near-optimal control profile. 
\subsection{Numerical Results for Case \#3}
\begin{figure}[!ht]
\centering
\includegraphics[width=3in]{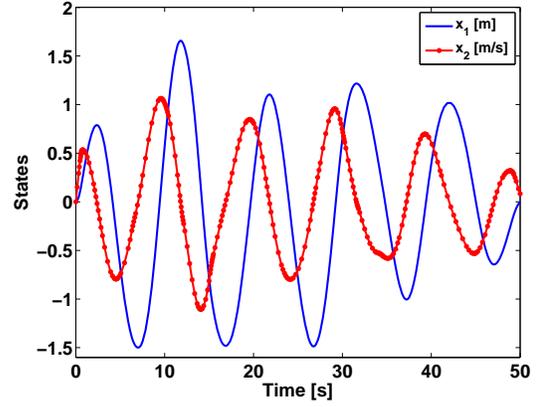}
\caption{Time histories of the states for Case 3.}
\label{fig:case3_states}
\end{figure}

\begin{figure}[!ht]
\centering
\includegraphics[width=3.4in]{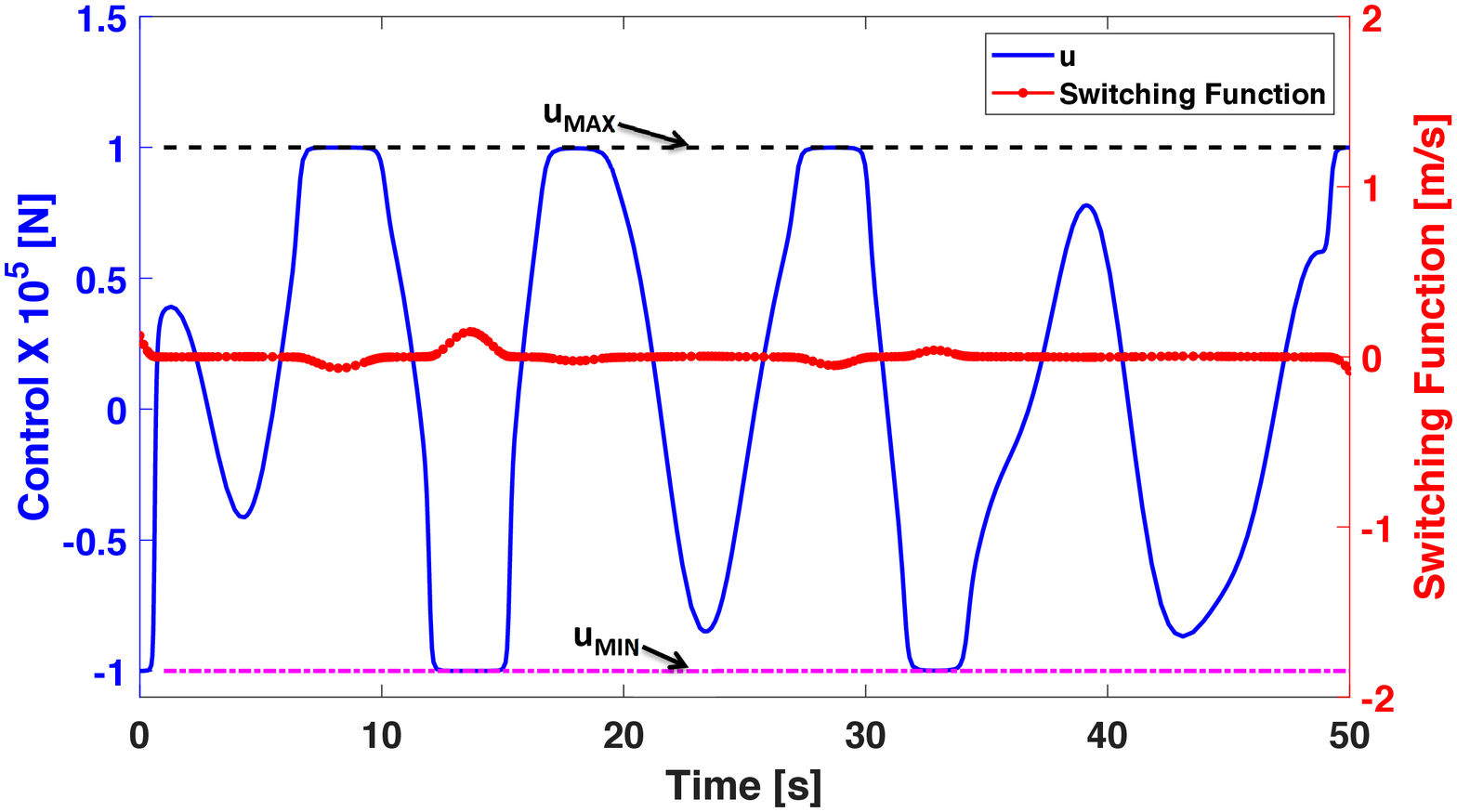}
\caption{Time history of control and switching function for Case 3.}
\label{fig:case3_control}
\end{figure}

\begin{figure}[!ht]
\centering
\includegraphics[width=2.8in]{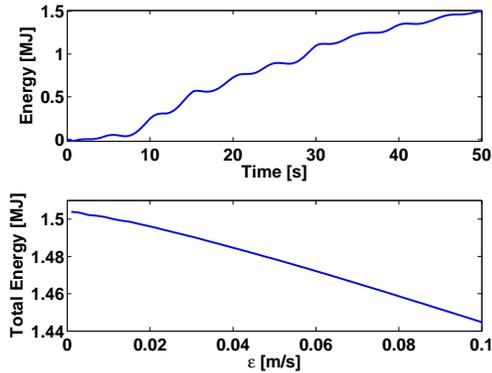}
\caption{Time history of harvested energy (upper plot) and change in the total harvested energy vs. $\epsilon_1$ for Case 3.}
\label{fig:case3_energy}
\end{figure}

The boundary conditions for this case are identical to those of case 2; Fig.~\ref{fig:case3_states} shows the states solutions. In this case, the bounds on the control are tighter as compared to cases 1 and 2. Note that only the results obtained using the ETRM are included for this case.
Figure \ref{fig:case3_control} shows the control profile for this case, which has a complicated structure comprising of the following sequence: B-S-B-S-B-S-B-S-B-S-B-S-B. Thus, the control solution has six singular arcs and seven bang arcs. Traditionally, based on Fig.~\ref{fig:case3_control}, one would have to solve a 14-point BVP using the indirect methods, which is more complicated than cases 1 and 2. However, the ETRM solves a simpler TPBVP for this case. The bounds on the controls are implemented in an automated and implicit manner using the ETRM. The switching function, shown in Fig.~\ref{fig:case3_control}, vanishes for the singular part, has positive values for the negative bang parts and negative values for the positive bang parts of the control solution. Thus, the switching function profile is in excellent agreement with the optimal control profile obtained for this case.  
Figure \ref{fig:case3_energy} shows the energy time history plot for this case. The values of $\epsilon$ for this case are decreased exactly like cases 1 and 2, leading to higher and more accurate values for the energy harvested from the waves. The energy harvested for case 3 is nearly twice the value of energy harvested in cases 1 and 2, which is due to the more energetic non-periodic excitation force in case 3. 
\subsection{Summary of Results}
The results for the three cases are summarized in Table~\ref{table:wecsummary}. The energy results for cases 1 and 2 demonstrate nearly 5\% and 6.5\% improvements, respectively, over the results from Ref. \cite{zou2017optimal}. The main difference between the results is caused by the terminal segment of the control, where the ETRM takes a bang form. The results in this study indicate that a final bang arc leads to a higher value for the absorbed energy and takes precedence over a singular arc.  
\begin{table}[!ht]
 \begin{center}
           \caption{Summary of results for different cases; solver relative and absolute tolerances are $1 \times 10^{-4}$ for all cases; $\bm{\epsilon}$ is $1 \times 10^{-3}$.}
 \begin{tabular}{p{2.65cm}cccc}
        \hline
        \hline
        \textbf{Attribute} & \textbf{Case 1} & \textbf{Case 2} & \textbf{Case 3}\\\hline
        \textbf{Energy (MJ)} & 0.8412 & 0.7599 & 1.5040\\
        \textbf{Energy (MJ) Ref. \cite{zou2017optimal}} & 0.7966 & 0.7166 & -\\
        \textbf{Computation Time (s)} & 152.84 & 157.95 & 135.48\\\hline
\end{tabular}
          \label{table:wecsummary}
 \end{center}
\end{table}
All computations were performed on a personal computer with a 2.5-GHz Intel i5 processor using MATLAB 2014b built-in BVP solver, \textit{bvp4c}.

\section{CONCLUSIONS}
\label{sec:conc}
In this study, the application of the Epsilon Trig Regularization Method (ETRM) is demonstrated to the problem of maximum-energy-absorption for a point-absorber WEC. The ETRM is a simple, efficient, and powerful method for dealing with OCPs with control constraints based on trigonometry. 
Using the ETRM, two trigonometric terms are added to the path cost of standard optimal control formulation and one of state equations to implement the regularization. 

We considered three scenarios for the WEC problem. The results indicate that high-quality and accurate solutions were obtained for these cases by using the ETRM as compared to the solutions obtained from the literature and a direct solver package based on pseudo-spectral methods (PSM). Singular control solutions obtained using the PSM involve many jitters, which are unrealistic to implement in a real world scenario. The results indicate that more wave energy can be harvested by using a combination of singular and bang control profiles (obtained using the ETRM) as compared to a purely singular control profile (as proposed in the literature). 

\bibliographystyle{IEEEtran} 
\bibliography{wectrig,RefET}


\end{document}